# Tie-respecting bootstrap methods for estimating distributions of sets and functions of eigenvalues

PETER HALL[1,2], YOUNG K. LEE[3], BYEONG U. PARK[4] and DEBASHIS PAUL[1]

[1]*Department of Statistics, University of California, One Shields Avenue, Davis, CA 95616, USA.*
*E-mail: debashis@wald.ucdavis.edu*

[2]*Department of Mathematics and Statistics, The University of Melbourne, Melbourne, VIC 3010, Australia. E-mail: halpstat@ms.unimelb.edu.au*

[3]*Department of Statistics, Kangwon National University, Chuncheon 200-701, Korea.*
*E-mail: youngklee@kangwon.ac.kr*

[4]*Department of Statistics, Seoul National University, Seoul 151-747, Korea.*
*E-mail: bupark@stats.snu.ac.kr*

Bootstrap methods are widely used for distribution estimation, although in some problems they are applicable only with difficulty. A case in point is that of estimating the distributions of eigenvalue estimators, or of functions of those estimators, when one or more of the true eigenvalues are tied. The $m$-out-of-$n$ bootstrap can be used to deal with problems of this general type, but it is very sensitive to the choice of $m$. In this paper we propose a new approach, where a tie diagnostic is used to determine the locations of ties, and parameter estimates are adjusted accordingly. Our tie diagnostic is governed by a probability level, $\beta$, which in principle is an analogue of $m$ in the $m$-out-of-$n$ bootstrap. However, the tie-respecting bootstrap (TRB) is remarkably robust against the choice of $\beta$. This makes the TRB significantly more attractive than the $m$-out-of-$n$ bootstrap, where the value of $m$ has substantial influence on the final result. The TRB can be used very generally; for example, to test hypotheses about, or construct confidence regions for, the proportion of variability explained by a set of principal components. It is suitable for both finite-dimensional data and functional data.

## 1. Introduction

Bootstrap methods can be particularly effective in distribution estimation, but typically only in cases where the distribution being estimated is asymptotically normal. Incon-sistency occurs in many settings, ranging from inference about extremes to problems involving goodness-of-fit testing.

Arguably the most commonly occurring difficulties of this type arise when estimating potentially tied eigenvalues. In this paper we suggest an adaptive, two-stage approach to tackling this problem, based on a new bootstrap algorithm. We show that a good bound for the distances between eigenvalues and their estimators, founded on an inequality borrowed from mathematical analysis, can be combined with the conventional bootstrap to give an effective statistical diagnostic for identifying places where ties occur. Armed with this information, a new, tie-respecting bootstrap algorithm can be employed to generate data that reflect the conclusion of the first bootstrap step.

Numerical and theoretical properties of the resulting tie-respecting bootstrap (TRB) are developed. Together they show that the method can be used reliably in a wide range of settings. Our theoretical contributions include a new representation for the limiting joint distribution of eigenvalue estimators, valid very generally – for example, in the functional data case and in both tied and untied eigenvalue settings.

A variety of diagnostics can be used in the first stage of the algorithm. The one on which we focus is "tuned" using a probability level, $\beta$; the TRB is remarkably robust against the choice of this quantity. This contrasts markedly with the $m$-out-of-$n$ bootstrap, which is particularly sensitive to the value of $m$. We demonstrate this point in a simulation study and by proving theoretically that the second stage of the TRB algorithm is largely unaffected by the nature of the diagnostic in the first stage; see the first paragraph of Section 4.2. We also show that a simple inequality provides a conservative way of accommodating the value chosen for $\beta$; see the last paragraph of Section 4.3.

The TRB is valid in conventional, finite-dimensional settings, where eigenvalues are defined in terms of matrices, and also in less standard problems involving functional data analysis. Since a non-expert reader may be unable to develop methodology in the functional data case, we introduce our methodology there. Its vector-case version is entirely analogous, and is treated briefly, but specifically, in Section 2.7. Particularly in the case of functional data, principal components analysis is a popular way of reducing dimension and the sizes of eigenvalues convey a great deal of information about the amount of variability that is captured by relatively low-dimensional approximations.

The problem of bootstrap-based inference for eigenvalues has a history that is not much younger than that of the bootstrap itself. As early as 1985, and in the vector case, Beran and Srivastava (1985) discussed consistency of bootstrap methods for confidence regions, noting that consistency fails in the event of eigenvalue ties. To overcome this difficulty they suggested avoiding the problem of computing confidence regions for individual eigenvalues, and constructing instead a simultaneous region for all the eigenvalues.

Alemayehu (1988) discussed techniques for finding approximate simultaneous confidence sets for functions of eigenvalues and eigenvectors. Beran (1988) developed a general bootstrap approach to constructing simultaneous confidence regions, and illustrated its application using the example of simultaneous regions for eigenvalues. In a general but parametric setting, Beran (1997) suggested selecting and adjusting parameter values of the distribution from which bootstrap data are drawn, in order to achieve consistency. Andrews (2000) developed a theory describing circumstances where bootstrap performance is compromised.

Other work on properties of resampling methods for inference about eigenvalues and eigenvectors includes that of Nagao (1988), who obtained limiting distributions of jackknife statistics associated with eigenvalue and eigenvector estimation; Eaton and Tyler (1991), who introduced techniques for deriving the asymptotic distributions of eigenvalue estimators, and illustrated them in the context of the bootstrap; Dümbgen (1993, 1995), who discussed bootstrap-based methods for confidence regions and hypothesis tests related to eigenvalues and eigenvectors; Zhang and Boos (1993), who introduced bootstrap tests of hypotheses about multivariate covariance structures; and Schott (2006), who suggested a test for equality of the smallest eigenvalues of a covariance matrix.

Hall *et al.* (1993) described estimation of the largest eigenvalue using the $m$-out-of-$n$ bootstrap in the case of ties. Bickel *et al.* (1997) and Bickel (2003) also discussed the $m$-out-of-$n$ bootstrap and, in particular, addressed its performance when used to estimate distributions that cannot be accessed using the standard bootstrap. However, when used in tied-eigenvalue problems, this technique is uncompetitive with the approach suggested here on two grounds: First, it requires empirical choice of $m$ for which a suitable algorithm does not seem to be available. Second, it produces distribution estimators that converge relatively slowly.

There is a large amount of literature, too, on principal component analysis for functional data. In that setting, methodology goes back at least to the work of Besse and Ramsay (1986), Ramsay and Dalzell (1991) and Rice and Silverman (1991). The literature is surveyed in greater detail by Ramsay and Silverman (2002, 2005). Relatively theoretical contributions to the functional-data case include those of Dauxois *et al.* (1982), Bosq (1989, 2000) and Besse (1992).

## 2. Methodology

### 2.1. Background: conventional estimators of eigenvalues and eigenvectors

Given a random sample $\mathcal{X} = \{X_1, \ldots, X_n\}$ from the distribution of a random function $X$, let $\widehat{K}(u,v)$ denote the conventional estimator of the covariance function, $K(u,v) = \operatorname{cov}\{X(u), X(v)\}$:

$$\widehat{K}(u,v) = \frac{1}{n} \sum_{i=1}^{n} \{X_i(u) - \bar{X}(u)\}\{X_i(v) - \bar{X}(v)\}, \tag{2.1}$$

where $\bar{X} = n^{-1} \sum_i X_i$. It will be assumed that the argument of $X$ is confined to a compact interval $\mathcal{I}$, say, and that $u$ and $v$ are also restricted to that region.

The eigenvalues $\theta_j$ and eigenvectors, or eigenfunctions, $\psi_j$, are arguably most clearly expressed in terms of the spectral decomposition of the linear operator of which the kernel is $K$:

$$K(u,v) = \sum_{j=1}^{\infty} \theta_j \psi_j(u) \psi_j(v). \tag{2.2}$$

Specifically, denoting the operator too by $K$, the operator is defined by $(K\psi)(u) = \int K(u,v)\psi(v)\,dv$, and in these terms, $K\psi_j = \theta_j \psi_j$. Here and below, unqualified integrals are taken over the interval $\mathcal{I}$.

The estimator $\widehat{K}$, at (2.1), admits an expansion analogous to that for $K$, at (2.2):

$$\widehat{K}(u,v) = \sum_{j=1}^{\infty} \hat{\theta}_j \hat{\psi}_j(u) \hat{\psi}_j(v). \tag{2.3}$$

In both (2.2) and (2.3) the eigenvalues are assumed to be ordered as decreasing sequences:

$$\theta_1 \geq \theta_2 \geq \cdots \geq 0, \qquad \hat{\theta}_1 \geq \hat{\theta}_2 \geq \cdots \geq 0. \tag{2.4}$$

The inherent positive semi-definiteness of a covariance function guarantees the non-negativity claimed in (2.4). The fact that the sample $\mathcal{X}$ contains only $n$ elements ensures that $\hat{\theta}_j$ vanishes for $j > n$. This in turn implies that the functions $\hat{\psi}_j$ in (2.3) are not determined for $j > n$.

Under mild continuity assumptions on the distribution of the stochastic process $X$, random fluctuations within the data set $\mathcal{X}$ guarantee that, even if there are ties among non-zero eigenvalues in the sequence $\theta_j$, these are not reflected among the empirical eigenvalues $\hat{\theta}_j$, with the result that $\hat{\theta}_1 > \cdots > \hat{\theta}_n$ with probability 1.

### 2.2. Principal components

The principal components of $X$ are the coefficients $\xi_j = \int (X - EX)\psi_j$, and lead to the Karhunen–Loève expansion,

$$X - \mathrm{E}(X) = \sum_{j=1}^{\infty} \xi_j \psi_j. \tag{2.5}$$

The definition of $\xi_j$ implies that those quantities are uncorrelated.

Analogously, the empirical principal components are defined by $\hat{\xi}_{ij} = \int (X_i - \bar{X})\hat{\psi}_j$, and lead to an empirical version of (2.5),

$$X_i - \bar{X} = \sum_{j=1}^{\infty} \hat{\xi}_{ij} \hat{\psi}_j. \tag{2.6}$$

Reflecting the properties $\mathrm{E}(\xi_j) = 0$ and $\mathrm{var}(\xi_j) = \theta_j$ enjoyed by the true principal components, we have for their empirical counterparts,

$$\sum_{i=1}^{n} \hat{\xi}_{ij} = 0, \qquad \frac{1}{n}\sum_{i=1}^{n} \hat{\xi}_{ij}^2 = \hat{\theta}_j, \tag{2.7}$$

for each $j$. Since only $n$ data curves $X_i$ are available, then $\hat{\xi}_{ij} = 0$ for $j \geq n+1$ and for each $i$.

### 2.3. A tie diagnostic

If there are ties among the $\theta_j$'s, then the corresponding values of $\hat\theta_j$ are generally $n^{-1/2}$ apart. Specifically, if $\theta_{p+1} = \cdots = \theta_{p+q}$, then, under regularity conditions, the differences $n^{1/2}(\hat\theta_{p+j} - \hat\theta_{p+k})$ have proper, non-degenerate limiting distributions for $1 \le j < k \le q$.

Herein lies the difficulty that conventional bootstrap methods have reflecting tied eigenvalues. In the standard bootstrap approximation to "reality", the $\hat\theta_j$'s represent the respective "true" eigenvalues $\theta_j$, and so the bootstrap incurs errors of size $n^{-1/2}$ when it is employed, explicitly or implicitly, to approximate the differences between identical eigenvalues. That is, in places where the eigenvalue differences should be zero, their values in the bootstrap world are of size $n^{-1/2}$. This is the same order as the difference between an eigenvalue estimator and the true eigenvalue; the extra term of size $n^{-1/2}$, representing a quantity that should really be zero, confounds the distribution-estimation problem. In consequence, the bootstrap estimator of the distribution of a tied eigenvalue is not consistent.

We suggest overcoming this problem by, first, using the data to estimate where the tied eigenvalues are, and subsequently replacing these empirically-determined ties by tied eigenvalue estimators. In principle, estimating the locations of ties requires us to have good estimators of the distributions of eigenvalue estimators, and that in turn demands knowledge of the locations of ties. We may break this circular argument by using a relatively robust method for constructing simultaneous confidence bands, such as that given below. Depending on the number of ties, and their locations in the eigenvalue sequence, our method can be improved by using a more sophisticated approach to constructing simultaneous bounds. However, the principle remains the same.

It is known that, with probability 1,

$$\sup_{j \ge 1} |\hat\theta_j - \theta_j| \le \|\!|\widehat K - K\|\!|, \tag{2.8}$$

where, for any bivariate function $L$, $\|\!| L \|\!|^2 = \iint L^2$. This property suggests that simultaneous confidence intervals for the $\theta_j$'s can be constructed by using a bootstrap procedure for estimating the distribution of $\|\!|\widehat K - K\|\!|$. To this end, let $X_1^\dagger, \ldots, X_n^\dagger$ denote a bootstrap resample drawn by sampling randomly, with replacement, from $\mathcal{X}$ in the conventional way, and put

$$\widehat K^\dagger(u,v) = \frac{1}{n} \sum_{i=1}^n \{X_i^\dagger(u) - \bar X^\dagger(u)\}\{X_i^\dagger(v) - \bar X^\dagger(v)\}, \tag{2.9}$$

where $\bar X^\dagger = n^{-1} \sum_i X_i^\dagger$. Choosing a probability level $\beta$, such as $\beta = 0.05$, take $\hat z_\beta$ to be the solution of the equation

$$P(\|\!|\widehat K^\dagger - \widehat K\|\!| \le \hat z_\beta | \mathcal{X}) = 1 - \beta.$$

Then, approximate and often slightly conservative simultaneous confidence bounds for $\theta_j$ are given by $\hat\theta_j \pm \hat z_\beta$ for each $j \ge 1$. Properties of this method were explored by Hall and Hosseini–Nasab (2006), who showed that the level of conservatism is usually slight.

If the confidence intervals $(\hat{\theta}_j - \hat{z}_\beta, \hat{\theta}_j + \hat{z}_\beta)$ and $(\hat{\theta}_{j+1} - \hat{z}_\beta, \hat{\theta}_{j+1} + \hat{z}_\beta)$ intersect, that is, if $\hat{\theta}_j - \hat{\theta}_{j+1} < 2\hat{z}_\beta$, then our tie diagnostic asserts that $\theta_j = \theta_{j+1}$. On the other hand, if $\hat{\theta}_j - \hat{\theta}_{j+1} \geq 2\hat{z}_\beta$, then our diagnostic states that $\theta_j > \theta_{j+1}$. These conclusions, which amount to the results of a sequence of simultaneous hypothesis tests, uniquely define a sequence of values of $\hat{p}_k$ and $\hat{q}_k$ for $k \geq 1$ starting with $\hat{p}_1 = 0$, where $\hat{p}_k + 1$ and $\hat{p}_k + \hat{q}_k$ define the end-points of the $k$th sequence of integers, $j$, for which the diagnostic asserts that the eigenvalues $\theta_j$ are equal to one another. Note that $\hat{q}_k = 1$ if the diagnostic suggests that there are no ties for $\theta_{\hat{p}_k+1}$. It follows that

$$\hat{p}_k + \hat{q}_k = \hat{p}_{k+1} \text{ for each } k \geq 1 \text{ and } 1 = \hat{p}_1 + 1 \leq \hat{p}_1 + \hat{q}_1 < \hat{p}_2 + 1 \leq$$
$$\hat{p}_2 + \hat{q}_2 < \hat{p}_3 + 1 \leq \hat{p}_3 + \hat{q}_3 < \cdots, \text{ with the sequence of inequalities} \quad (2.10)$$
$$\text{ending when we find a value of } \nu \text{ for which } \hat{q}_\nu = \infty.$$

One may define an alternative tie diagnostic based on the bootstrap distribution of $\sup_{j \geq 1} |\hat{\theta}_j - \theta_j|$. Let $\hat{\theta}_j^\dagger$ denote the eigenvalues of the covariance operator $\widehat{K}^\dagger$, and $\tilde{z}_\beta$ be the solution of the equation

$$P\left(\sup_{j \geq 1} |\hat{\theta}_j^\dagger - \hat{\theta}_j| \leq \tilde{z}_\beta \Big| \mathcal{X}\right) = 1 - \beta.$$

The corresponding tie diagnostic states that $\theta_j = \theta_{j+1}$ if $\hat{\theta}_j - \hat{\theta}_{j+1} < 2\tilde{z}_\beta$; $\theta_j > \theta_{j+1}$, otherwise.

### 2.4. Adjusting eigenvalue estimators to reflect ties

Suppose that an empirical tie diagnostic, such as the bootstrap-based method discussed in Section 2.3, suggests that ties occur among $\theta_{\hat{p}_k+1}, \ldots, \theta_{\hat{p}_k+\hat{q}_k}$ for $k \geq 1$, where $\hat{p}_k$ and $\hat{q}_k$ satisfy (2.10). Then we modify the set of eigenvalue estimators as follows:

(a) For j in the range $\hat{p}_k + 1 \leq j \leq \hat{p}_k + \hat{q}_k$, we replace $\hat{\theta}_j$ by the average, $\tilde{\theta}_j$ say, of the values of $\hat{\theta}_{\hat{p}_k+1}, \ldots, \hat{\theta}_{\hat{p}_k+\hat{q}_k}$, provided $1 \leq k \leq \nu - 1$. (b) We replace $\hat{\theta}_j$ by the average value, $\tilde{\theta}_j$, of $\hat{\theta}_{\hat{p}_\nu+1}, \ldots, \hat{\theta}_n$ if $\hat{p}_\nu + 1 \leq j \leq n$; and we leave the value of $\hat{\theta}_j$ unchanged at zero if $j > n$, but relabel it $\tilde{\theta}_j$. (2.11)

In addition to producing ties in the estimated eigenvalue sequence when the tie diagnostic says they should be there, part (a) of the algorithm at (2.11) identifies "ties of order one", that is, instances where $\hat{q}_k = 1$ and the corresponding value of $\hat{\theta}_j$ (with $j = \hat{p}_k + 1 = \hat{p}_k + \hat{q}_k$) is equal to $\tilde{\theta}_j$. In part (b) we could have replaced $\hat{\theta}_j$ by zero whenever $j \geq \hat{p}_\nu + 1$, but that would have altered the total value, $\sum_j \hat{\theta}_j$, of the estimated eigenvalues. This quantity estimates the total variability of the random function $X$, and is of statistical importance in its own right, without regard to individual eigenvalue estimators. Therefore, we would prefer to leave it unchanged.

Implementing the algorithm at (2.11), we generate a new sequence $\tilde{\theta}_1 \geq \tilde{\theta}_2 \geq \cdots$ of eigenvalue estimators.

## 2.5. Correcting the empirical principal components

It can be deduced from the second part of (2.7) that, after adjusting the values of the eigenvalue estimators $\hat{\theta}_j$ to reflect our assessment of where ties lie among the true eigenvalues, we should also rescale the empirical principal components before resampling. Thus, we are led to work with

$$\tilde{\xi}_{ij} = (\tilde{\theta}_j/\hat{\theta}_j)^{1/2}\hat{\xi}_{ij} \qquad \text{for } 1 \leq i \leq n \text{ and } 1 \leq j \leq n. \tag{2.12}$$

These quantities satisfy

$$\sum_{i=1}^{n} \tilde{\xi}_{ij} = 0, \qquad \frac{1}{n}\sum_{i=1}^{n} \tilde{\xi}_{ij}^2 = \tilde{\theta}_j.$$

Reflecting the uncorrected case, we define $\tilde{\xi}_{ij} = 0$ for $j \geq n+1$ and for each $i$.

Once we have computed these "corrected" principal components, we resample their values in much the same way that we would resample from a set of residuals in a regression problem. This is unusual in bootstrap algorithms for functional data analysis; usually the raw data are resampled. However, if one were to work instead with the conventional principal components $\hat{\xi}_{ij}$, rather than their corrected counterparts $\tilde{\xi}_{ij}$, then resampling, using the principles outlined below, would be equivalent to resampling from the raw data. Although the eigenfunctions can change signs and so can the principal components, this does not cause any difficulty since the resampled vectors of the principal components are multiplied to the empirical eigenfunctions $\hat{\psi}_j$ to produce a bootstrap sample.

## 2.6. A tie-respecting bootstrap algorithm

Let $\tilde{\xi}_i = (\tilde{\xi}_{i1}, \tilde{\xi}_{i2}, \ldots)$ denote the vector of corrected empirical principal components, the latter defined at (2.12). Conditional on the data $\mathcal{X}$, draw a resample $\tilde{\xi}_1^*, \ldots, \tilde{\xi}_n^*$ by sampling randomly, with replacement, from the collection $\tilde{\xi}_1, \ldots, \tilde{\xi}_n$. Put

$$X_i^* = \bar{X} + \sum_{j=1}^{\infty} \tilde{\xi}_{ij}^* \hat{\psi}_j,$$

this being a bootstrap version of the conventional Karhunen–Loève expansion at (2.5). We shall construct percentile-method bootstrap confidence regions for eigenvalues, using the resampled data $X_i^*$.

The bootstrap version of $\widehat{K}$, at (2.1), is given by

$$\widehat{K}^*(u,v) = \frac{1}{n}\sum_{i=1}^{n}\{X_i^*(u) - \bar{X}^*(u)\}\{X_i^*(v) - \bar{X}^*(v)\},$$

where $\bar{X}^* = n^{-1} \sum_i X_i^*$. (Thus, $\widehat{K}^*$ is the TRB form of $\widehat{K}$, whereas $\widehat{K}^\dagger$, at (2.9), is the conventional bootstrap form.) A bootstrap analogue of the spectral expansion (2.3) is

$$\widehat{K}^*(u,v) = \sum_{j=1}^{\infty} \hat{\theta}_j^* \hat{\psi}_j^*(u) \hat{\psi}_j^*(v).$$

The algorithm (2.11) can be applied directly to the bootstrap quantities $\hat{\theta}_k^*$, as at (2.13) below, just as it was earlier to the non-bootstrap empirical values. In the bootstrap case we employ, for simplicity, the same values $\hat{p}_k$ and $\hat{q}_k$ determined by a procedure we gave in Section 2.3.

(a) For $j$ in the range $\hat{p}_k + 1 \leq j \leq \hat{p}_k + \hat{q}_k$, we replace $\hat{\theta}_j^*$ by the average, $\tilde{\theta}_j^*$ say, of the values of $\hat{\theta}_{\hat{p}_k+1}^*, \ldots, \hat{\theta}_{\hat{p}_k+\hat{q}_k}^*$, provided $1 \leq k \leq \nu - 1$. (b) We replace $\hat{\theta}_j^*$ by the average value, $\tilde{\theta}_j^*$, of $\hat{\theta}_{\hat{p}_\nu+1}^*, \ldots, \hat{\theta}_n^*$ if $\hat{p}_\nu + 1 \leq j \leq n$, and we leave the value of $\hat{\theta}_j^*$ unchanged at zero if $j > n$, but relabel it $\tilde{\theta}_j^*$. (2.13)

In (2.13) it might be more natural to replace $\hat{p}_j$ and $\hat{q}_j$ by their respective bootstrap versions, $\hat{p}_j^*$ and $\hat{q}_j^*$. However, the error incurred by not doing this will generally be small.

To compute percentile-bootstrap confidence regions for $\theta_j$, we first solve for $\hat{x}_{j\alpha}$, as nearly as possible, the equation

$$P(\tilde{\theta}_j^* \leq \tilde{\theta}_j + \hat{x}_{j\alpha} | \mathcal{X}) = \alpha.$$

Here, $\alpha \in (0,1)$ represents a probability. Tie-respecting, one- and two-sided percentile-method confidence regions, each with nominal coverage $1 - \alpha$, are given by

$$(\tilde{\theta}_j - \hat{x}_{j,1-\alpha}, \infty), \qquad (-\infty, \tilde{\theta}_j - \hat{x}_{j\alpha}), \qquad (\tilde{\theta}_j - \hat{x}_{j,1-(\alpha/2)}, \tilde{\theta}_j - \hat{x}_{j,\alpha/2}). \qquad (2.14)$$

A simultaneous confidence region for a general, finite sequence of eigenvalues $\theta_{j_1}, \ldots, \theta_{j_k}$ can be constructed analogously. This region will have asymptotically correct coverage, even if some of the eigenvalues $\theta_{j_\ell}$ are tied with one another or with other eigenvalues not included in the sequence.

TRB confidence regions for functions of the eigenvalues $\theta_1, \theta_2, \ldots$ can be constructed using the same procedure. We illustrate below in the case of a confidence region for the ratio,

$$\rho = \left( \sum_{j=1}^{k} \theta_j \right) \bigg/ \left( \sum_{j=1}^{\infty} \theta_j \right),$$

which represents the proportion of the variability of the random function $X$ that is explained by the first $k$ principal components.

Define the tie-respecting estimator $\tilde{\rho}$ of $\rho$, and its bootstrap version $\tilde{\rho}^*$, by

$$\tilde{\rho} = \left( \sum_{j=1}^{k} \tilde{\theta}_j \right) \bigg/ \left( \sum_{j=1}^{\infty} \tilde{\theta}_j \right), \qquad \tilde{\rho}^* = \left( \sum_{j=1}^{k} \tilde{\theta}_j^* \right) \bigg/ \left( \sum_{j=1}^{\infty} \tilde{\theta}_j^* \right).$$

Solve for $\hat{y}_\alpha$, as nearly as possible, the equation

$$P(\tilde{\rho}^* \leq \tilde{\rho} + \hat{y}_\alpha | \mathcal{X}) = \alpha.$$

Analogously to (2.14), nominal $(1-\alpha)$-level, one- and two-sided percentile-method confidence regions for $\rho$ are given by

$$(\tilde{\rho} - \hat{y}_{1-\alpha}, \infty), \qquad (-\infty, \tilde{\rho} - \hat{y}_\alpha), \qquad (\tilde{\rho} - \hat{y}_{1-(\alpha/2)}, \tilde{\rho} - \hat{y}_{\alpha/2}). \tag{2.15}$$

In practice, the quantile $\hat{z}_\beta$ as defined in Section 3 is approximated by the corresponding quantile of the "empirical" distribution of $\|\widehat{K}^\dagger - \widehat{K}\|$ obtained from a finite number of bootstrap resamples. Different resamples yield different approximations of $\hat{z}_\beta$, and thus produce differing numbers of distinct eigenvalues when using the tie determination rule. This may be a source of additional variability in our TRB procedure and may influence the coverage performance of the TRB confidence intervals. The additional variability in $\hat{p}_k$ and $\hat{q}_k$ may be significant depending on the value of $\beta$ and spacings of eigenvalues. If this were of concern then it could be dealt with by using Breiman's (1996) bagging method, even to the extent of deleting simulations that did not accord with the majority assessment of the number of eigenvalues. We shall not explore this approach, however. One promising aspect of the present TRB method, as found in the numerical study presented in Section 3, is that the method is fairly robust against the choice of $\beta$ and works pretty well in various settings of spacings of eigenvalues.

## 2.7. Adaptation to matrix setting

If the data $X_i$ are random $p$-vectors, rather than random functions, then $K$ and $\widehat{K}$ should be interpreted as $p \times p$ matrices, with $\widehat{K}$ given by:

$$\widehat{K} = \frac{1}{n} \sum_{i=1}^{n} (X_i - \bar{X})(X_i - \bar{X})^{\mathrm{T}}.$$

The eigenvectors $\psi_j$ are now $p$-vectors, and infinite expansions, for example, the spectral decompositions at (2.2) and (2.3) are now of length only $p$:

$$K = \sum_{j=1}^{p} \theta_j \psi_j \psi_j^{\mathrm{T}}, \qquad \widehat{K} = \sum_{j=1}^{p} \hat{\theta}_j \hat{\psi}_j \hat{\psi}_j^{\mathrm{T}}.$$

Principal components are defined by vector multiplication rather than integration. In particular, $\hat{\xi}_{ij} = (X_i - \bar{X})^{\mathrm{T}} \hat{\psi}_j$. Once these reinterpretations are made, the account of the TRB algorithm in Section 2.5 is applicable to the case of vector-valued data.

The $L_2$ norm for bivariate functions $L$, defined in Section 2.6 by $\|L\|^2 = \iint L^2$, is given in the matrix case, where $L = (\ell_{i_1 i_2})$, say, by $\|L\|^2 = \sum_{i_1} \sum_{i_2} \ell_{i_1 i_2}^2$. With this reinterpretation of notation, (2.8) holds in its original form, and the discussion in Section 2.6 of both the tie diagnostic and its bootstrap-based implementation is valid in the vector case.

## 3. Numerical properties

Here we assess the finite sample performance of the proposed TRB algorithm described in Section 2.6. The sample functions $X_i$, $i = 1, \ldots, n$ were generated from the model

$$X(u) = \sum_{j=1}^{\infty} \xi_j \psi_j(u), \qquad u \in \mathcal{I},$$

where $\mathcal{I} = [-1, 1]$, $\xi_j$ were independently distributed as $N(0, \theta_j)$, and $\psi_j(u) = \sqrt{2}\cos(j\pi u)$. Recall that $\theta_j$ are the eigenvalues of the covariance operator $K$ defined in Section 2.1. We set $\theta_j = \{500 + 100(j-4)\}^{-1}$ for $4 \leq j \leq n$ and $\theta_j = 0$ for $j > n$ and considered the following three models for the values of $\theta_1, \theta_2$ and $\theta_3$:

(1) $\theta_1 = \theta_2 = \theta_3 = 1$,
(2) $\theta_1 = 1.6, \theta_2 = \theta_3 = 0.7$,
(3) $\theta_1 = 1.6, \theta_2 = 1, \theta_3 = 0.4$.

In these models the first three principal components explain most of the variance of $X$. We considered two sample sizes, $n = 100$ and 400. The proportion of the variability explained by the first three principal components equals 99.0% when $n = 100$ and 98.5% when $n = 400$.

To compute the estimators $\hat{\theta}_j$, at (2.3), we discretized the sample curves $X_i$. In particular, we took $J(\leq n)$ equi-spaced points $u_1, \ldots, u_J$ on $\mathcal{I}$, and performed a singular-value decomposition on the $n \times J$ matrix $Z = (X_i(u_j) - \bar{X}(u_j))$. This gave the estimators $\hat{\vartheta}_1 \geq \cdots \geq \hat{\vartheta}_J \geq 0$ of $Z^T Z$, and their corresponding eigenvectors $\hat{\phi}_j = (\hat{\phi}_{j1}, \ldots, \hat{\phi}_{jJ})^T$. The eigenvalues $\hat{\theta}_j$ and eigenfunctions $\hat{\psi}_j$ were then obtained from the formulae $\hat{\theta}_j = \hat{\vartheta}_j/(nJ)$ and $\hat{\psi}_j(s_k) = \sqrt{J}\hat{\phi}_{jk}$. The principal components $\hat{\xi}_{ij} = \int (X_i - \bar{X})\hat{\psi}_j$ were approximated by discretizing the integrals. In our numerical experiments we took $J = 30$ when $n = 100$ and $J = 100$ for $n = 400$. In respect to these calculation details, and others given below, we used the same settings in the "bootstrap world" as in non-bootstrap cases.

To assess performance we considered the coverage probabilities of the two-sided confidence intervals for $\theta_j$, at (2.14), and for the ratios $\rho_j = \sum_{1 \leq k \leq j} \theta_k / \sum_{k \geq 1} \theta_k$, at (2.15). We investigated the performance of the two tie diagnostics described in Section 2.3, one based on the bootstrap distribution of $\|\widehat{K} - K\|$ and the other on the bootstrap distribution of $\sup_{j \geq 1} |\hat{\theta}_j - \theta_j|$. In addition to the coverage probabilities, we evaluated how well these tie diagnostics identified the tied eigenvalues. Note that for model (1), $p_k = k + 1$ for $k \geq 2$. For model (2), $p_2 = 1$ and $p_k = k$ for $k \geq 3$, and for model (3), $p_k = k - 1$ for $k \geq 2$. In all cases, $p_1 = \hat{p}_1 = 0$ by definition. We computed $P(\hat{p}_2 = 3)$ for model (1), $P(\hat{p}_2 = 1, \hat{p}_3 = 3)$ for model (2) and $P(\hat{p}_2 = 1, \hat{p}_3 = 2, \hat{p}_4 = 3)$ for model (3). These are the probabilities of identifying correctly the ties among $\theta_j$ for $1 \leq j \leq 4$. Those four terms are critical to the success of the tie diagnostics method, since $\theta_j$, for $j \geq 4$, makes a negligible contribution to the total variation of $X$.

We also compared our method with $m$-out-of-$n$ bootstrap algorithms. The latter are based on bootstrap resamples of size $m(\leq n)$, drawn by sampling randomly with replacement from the sample $\mathcal{X} = \{X_1, \ldots, X_n\}$ in the conventional way. This approach ignores ties among the true eigenvalues. Let $\mathcal{X}_m^\dagger = \{X_1^\dagger, \ldots, X_m^\dagger\}$ be the $m$-out-of-$n$ bootstrap resample and $\widehat{K}_m^\dagger$ be a version of $\widehat{K}^\dagger$, defined in Section 2.3, that is based on $\mathcal{X}_m^\dagger$, that is,

$$\widehat{K}_m^\dagger(u,v) = \frac{1}{m} \sum_{i=1}^m \{X_i^\dagger(u) - \bar{X}_m^\dagger(u)\}\{X_i^\dagger(v) - \bar{X}_m^\dagger(v)\},$$

where $\bar{X}_m^\dagger = m^{-1} \sum_{i \leq m} X_i^\dagger$. Denote by $\hat{\theta}_{j,m}^\dagger$ the eigenvalues of the covariance operator $\widehat{K}_m^\dagger$ and let $\hat{x}_{j,\alpha,m}^\dagger$ be the solution of the equation

$$P(\hat{\theta}_{j,m}^\dagger \leq \hat{\theta}_j + \hat{x}_{j,\alpha,m}^\dagger | \mathcal{X}) = \alpha.$$

Then $m$-out-of-$n$ percentile-bootstrap confidence regions can be obtained by replacing $\tilde{\theta}_j$ and $\hat{x}_{j\alpha}$ in Section 2.6 by $\hat{\theta}_j$ and $(m/n)^{1/2}\hat{x}_{j,\alpha,m}^\dagger$, respectively. The normalization $(m/n)^{1/2}$ here derives from the fact that the conditional distribution of $m^{1/2}(\hat{\theta}_{j,m}^\dagger - \hat{\theta}_j)$, given $\mathcal{X}$, approximates the distribution of $n^{1/2}(\hat{\theta}_j - \theta_j)$. Results for $n = 400$ are reported in Tables 1–3, which contain coverage probabilities of the confidence regions at the nominal level $1 - \alpha = 0.9$.

**Table 1.** Coverage probabilities of confidence intervals at nominal level 0.9 for the model (1)

| Method | | Coverage probabilities for | | | | | |
|---|---|---|---|---|---|---|---|
| | | $\theta_1$ | $\theta_2$ | $\theta_3$ | $\rho_1$ | $\rho_2$ | $\tau$ |
| $m$-out-of-$n$ | $m/n = 1$ | 0.818 | 0.860 | 0.598 | 0.670 | 0.578 | |
| bootstrap | $m/n = 3/4$ | 0.830 | 0.848 | 0.620 | 0.702 | 0.594 | |
| | $m/n = 1/2$ | 0.842 | 0.846 | 0.620 | 0.722 | 0.602 | |
| | $m/n = 1/4$ | 0.876 | 0.854 | 0.644 | 0.772 | 0.636 | |
| | $m/n = 1/8$ | 0.892 | 0.844 | 0.644 | 0.792 | 0.670 | |
| Tie-respecting | $\beta = 0.1$ | 0.902 | 0.902 | 0.902 | 0.892 | 0.892 | 1.000 |
| bootstrap | 0.3 | 0.902 | 0.902 | 0.902 | 0.892 | 0.892 | 1.000 |
| based on TD1 | 0.5 | 0.902 | 0.902 | 0.902 | 0.892 | 0.892 | 1.000 |
| | 0.7 | 0.902 | 0.902 | 0.898 | 0.886 | 0.886 | 0.992 |
| | 0.9 | 0.866 | 0.866 | 0.848 | 0.812 | 0.812 | 0.910 |
| Tie-respecting | $\beta = 0.1$ | 0.902 | 0.902 | 0.900 | 0.888 | 0.888 | 0.996 |
| bootstrap | 0.3 | 0.870 | 0.870 | 0.864 | 0.838 | 0.838 | 0.938 |
| based on TD2 | 0.5 | 0.812 | 0.820 | 0.762 | 0.668 | 0.668 | 0.742 |
| | 0.7 | 0.764 | 0.806 | 0.642 | 0.448 | 0.448 | 0.448 |
| | 0.9 | 0.786 | 0.852 | 0.576 | 0.532 | 0.468 | 0.138 |

Note. Based on 500 pseudo-samples of size $n = 400$. TD1 stands for the tie-diagnostic based on the bootstrap estimate of the distribution of $\|\widehat{K} - K\|$, and TD2 for the diagnostic based on the bootstrap estimate of the distribution of $\sup_{j \geq 1} |\hat{\theta}_j - \theta_j|$. The numbers in the rightmost column are the values of $\tau = P(\hat{p}_2 = 3)$, the probability of identifying correctly the ties among $\theta_j$ for $1 \leq j \leq 4$ in this case.

**Table 2.** Coverage probabilities of confidence intervals at nominal level 0.9 for the model (2)

| Method | | Coverage probabilities for | | | | | |
|---|---|---|---|---|---|---|---|
| | | $\theta_1$ | $\theta_2$ | $\theta_3$ | $\rho_1$ | $\rho_2$ | $\tau$ |
| $m$-out-of-$n$ | $m/n = 1$ | 0.864 | 0.888 | 0.756 | 0.882 | 0.786 | |
| bootstrap | $m/n = 3/4$ | 0.872 | 0.888 | 0.770 | 0.890 | 0.776 | |
| | $m/n = 1/2$ | 0.862 | 0.876 | 0.780 | 0.882 | 0.762 | |
| | $m/n = 1/4$ | 0.852 | 0.874 | 0.778 | 0.870 | 0.770 | |
| | $m/n = 1/8$ | 0.840 | 0.876 | 0.764 | 0.850 | 0.758 | |
| Tie-respecting | $\beta = 0.1$ | 0.866 | 0.880 | 0.880 | 0.884 | 0.878 | 0.994 |
| bootstrap | 0.3 | 0.866 | 0.886 | 0.886 | 0.884 | 0.880 | 1.000 |
| based on TD1 | 0.5 | 0.866 | 0.886 | 0.886 | 0.884 | 0.880 | 1.000 |
| | 0.7 | 0.866 | 0.886 | 0.886 | 0.884 | 0.880 | 1.000 |
| | 0.9 | 0.866 | 0.878 | 0.874 | 0.884 | 0.870 | 0.986 |
| Tie-respecting | $\beta = 0.1$ | 0.866 | 0.886 | 0.886 | 0.884 | 0.880 | 1.000 |
| bootstrap | 0.3 | 0.866 | 0.886 | 0.882 | 0.884 | 0.876 | 0.996 |
| based on TD2 | 0.5 | 0.864 | 0.872 | 0.870 | 0.884 | 0.860 | 0.954 |
| | 0.7 | 0.864 | 0.846 | 0.794 | 0.884 | 0.800 | 0.810 |
| | 0.9 | 0.864 | 0.860 | 0.744 | 0.884 | 0.776 | 0.456 |

Note. Based on 500 pseudo-samples of size $n = 400$. The meanings of TD1 and TD2 are the same as in Table 1. The numbers in the rightmost column are the values of $\tau = P(\hat{p}_2 = 1, \hat{p}_3 = 3)$, the probability of identifying correctly the ties among $\theta_j$ for $1 \leq j \leq 4$ in this case.

Section 2.3 described two tie diagnostics. In the discussion below and in the tables, we refer to the tie diagnostic based on the bootstrap distribution of $\|\widehat{K} - K\|$ as TD1, and the other as TD2. Let $\tau$ denote the probability of identifying correctly the ties among $\theta_j$ for $1 \leq j \leq 4$. The results in the tables show that the performance of the $m$-out-of-$n$ bootstrap is sensitive to the choice of $m$. When ties are present the subsampling scheme improves the coverage probability of the conventional bootstrap ($m = n$), but it suffers from the need to choose $m$ carefully. Furthermore, in the case of ties it has poor coverage accuracy even when $m$ is chosen to give the best performance.

In contrast, the TRB method with TD1 diagnostic is quite successful over the entire range of $\beta$ in these cases; that is, it is robust against the choice of $\beta$. Moreover, it has greater coverage accuracy than the $m$-out-of-$n$ bootstrap at optimal values of $\beta$ and $m$, respectively. When there are no ties among eigenvalues, the conventional bootstrap works well (as expected), but here the TRB methods perform well, too, exhibiting essentially the same coverage probabilities as the conventional bootstrap method for $\beta > 0.1$ in the case of TD2, and for $\beta > 0.3$ in the case of TD1.

Looking at the probability $\tau$ one finds that TD1 picks up the tied eigenvalues correctly, with very high probability, for all values of $\beta$, while TD2 works better than TD1 for small values of $\beta$ when there is no tie. These properties are directly translated to the performance of the methods in terms of coverage probability. Recall that $\sup_{j \geq 1} |\hat{\theta}_j - \theta_j| \leq$

**Table 3.** Coverage probabilities of confidence intervals at nominal level 0.9 for the model (3)

|  |  | Coverage probabilities for | | | | | |
|---|---|---|---|---|---|---|---|
| Method |  | $\theta_1$ | $\theta_2$ | $\theta_3$ | $\rho_1$ | $\rho_2$ | $\tau$ |
| $m$-out-of-$n$ | $m/n = 1$ | 0.860 | 0.884 | 0.876 | 0.878 | 0.864 | |
| bootstrap | $m/n = 3/4$ | 0.872 | 0.872 | 0.886 | 0.886 | 0.854 | |
|  | $m/n = 1/2$ | 0.866 | 0.868 | 0.884 | 0.866 | 0.854 | |
|  | $m/n = 1/4$ | 0.838 | 0.854 | 0.880 | 0.826 | 0.852 | |
|  | $m/n = 1/8$ | 0.814 | 0.828 | 0.870 | 0.774 | 0.830 | |
| Tie-respecting | $\beta = 0.1$ | 0.742 | 0.684 | 0.006 | 0.736 | 0.618 | 0.004 |
| bootstrap | 0.3 | 0.846 | 0.858 | 0.526 | 0.868 | 0.768 | 0.544 |
| based on TD1 | 0.5 | 0.858 | 0.878 | 0.866 | 0.876 | 0.860 | 0.974 |
|  | 0.7 | 0.860 | 0.884 | 0.876 | 0.878 | 0.862 | 1.000 |
|  | 0.9 | 0.860 | 0.884 | 0.876 | 0.878 | 0.862 | 1.000 |
| Tie-respecting | $\beta = 0.1$ | 0.856 | 0.868 | 0.584 | 0.872 | 0.780 | 0.634 |
| bootstrap | 0.3 | 0.860 | 0.884 | 0.876 | 0.878 | 0.862 | 0.998 |
| based on TD2 | 0.5 | 0.860 | 0.884 | 0.876 | 0.878 | 0.862 | 1.000 |
|  | 0.7 | 0.860 | 0.884 | 0.876 | 0.878 | 0.862 | 1.000 |
|  | 0.9 | 0.860 | 0.884 | 0.876 | 0.878 | 0.862 | 1.000 |

Note. Based on 500 pseudo-samples of size $n = 400$. The meanings of TD1 and TD2 are the same as in Table 1. The numbers in the rightmost column are the values of $\tau = P(\hat{p}_2 = 1, \hat{p}_3 = 2, \hat{p}_4 = 3)$, the probability of identifying correctly the ties among $\theta_j$ for $1 \leq j \leq 4$ in this case.

$\|\widehat{K} - K\|$ at (2.8). This implies that TD1 tends to accept $\theta_j = \theta_{j+1}$ more often than TD2, for fixed $\beta$.

In the case $n = 100$, for which the results are not reported in the tables, we found that the TRB does not work as well as when $n = 400$. Nevertheless, there exists a range of $\beta$ for each tie diagnostic such that the resulting TRB does better than the conventional bootstrap. One interesting point here is that the TRB recovers quite fast as the sample size increases to $n = 400$. This is true for all cases in the models (1) and (2) where there is a tie. For the model (3), one finds that, even for sample size $n = 400$, the coverage probabilities are still far away from their nominal value in a few cases. By increasing the sample size further to $n = 1000$, we found that the coverage probabilities in these cases improve greatly. For example, in the case where TD1 with $\beta = 0.1$ is employed, the coverage probabilities for $\theta_1, \theta_2, \theta_3, \rho_1, \rho_2$ equal $0.870, 0.872, 0.896, 0.874, 0.892$, respectively. One lesson learned from this observation is that the coverage error of the TRB converges to zero quite fast as $n$ increases, and the convergence depends only a little on the choice of $\beta$.

To see the performance of the TRB algorithm and of the $m$-out-of-$n$ bootstrap method as spacings of eigenvalues change on a continuous scale, we calculated the coverage probabilities of the two-sided confidence intervals for $\theta_j, 1 \leq j \leq 3$, when

$$(\theta_1, \theta_2, \theta_3) = (1 + s, 1, 1 - s), \qquad 0 \leq s \leq 0.5,$$

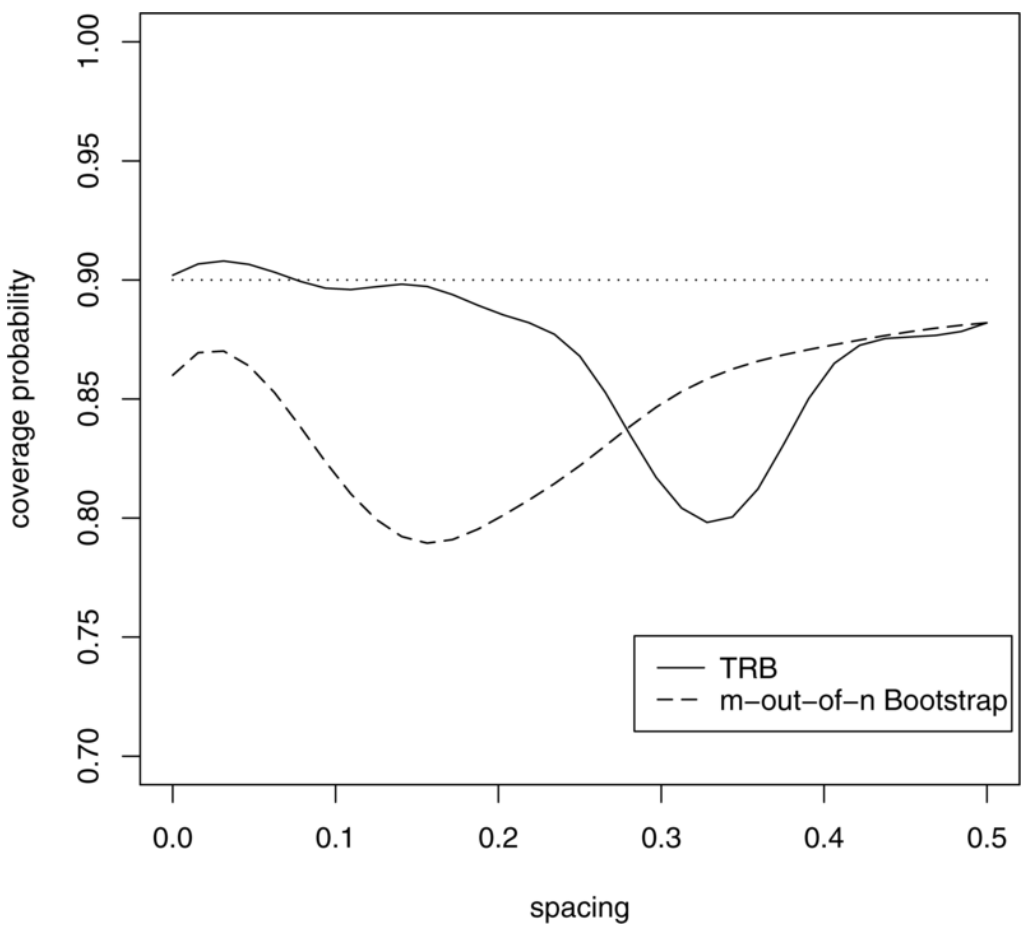

**Figure 1.** Coverage probabilities for $\theta_2 = 1$ as functions of the spacing $s = \theta_2 - \theta_1 = \theta_3 - \theta_2$ at nominal level 0.9, based on 500 pseudo-samples of size $n = 400$. The solid curve corresponds to TD1 with $\beta$ that gives the best coverage performance, and the dashed is for the $m$-out-of-$n$ with the optimal sampling fraction $m/n$.

and $\theta_j$ for $j \geq 4$ were the same as those in the three models (1)–(3). Figure 1 depicts the coverage probabilities of TD1 and the $m$-out-of-$n$ bootstrap method as functions of $s$, for which $\beta$ and $m/n$, respectively, were chosen to give the best coverage performance. These results suggest that the proposed TRB method is also robust against spacings of eigenvalues.

Although not reported here, we also implemented bootstrap calibration for each method for enhancing coverage accuracy of confidence regions. We found that the double bootstrap has little effect on coverage accuracy in the settings considered here. It produces confidence regions that have nearly the same coverage error as their non-calibrated versions when the latter do not do well. Calibration slightly improves coverage error, or makes it worse, when the coverage accuracy is good.

The main conclusions from the numerical study are that, in the presence of ties, the TRB works well for moderate sample sizes; that TD1 correctly identifies tied eigenvalues with high probability, and thus gives good coverage performance for a wide range of values of $\beta$ and in various settings of spacings of eigenvalues; and that, in the case of ties, this approach outperforms the $m$-out-of-$n$ bootstrap, in terms of coverage accuracy and robustness against choice of tuning parameter. We repeated our experiments for nominal level $1 - \alpha = 0.95$, and found that the lessons there were essentially the same.

## 4. Theoretical properties

### 4.1. Distributions of conventional estimators of eigenvalues

We first describe spacings among eigenvalues using a model that reflects the empirical description at (2.10):

$q_k \geq 1$ and $p_k + q_k = p_{k+1}$ for each $k$, and $1 = p_1 + 1 \leq p_1 + q_1 < p_2 + 1 \leq p_2 + q_2 < p_3 + 1 \leq p_3 + q_3 < \ldots$, with the sequence of inequalities either continuing ad infinitum, in which case we define $\nu = \infty$, or ending with a finite value of $\nu$ for which $q_\nu = \infty$; $\theta_{j_1} = \theta_{j_2}$ for $p_k + 1 \leq j_1, j_2 \leq p_k + q_k$ and for each $k$; and $\theta_{j_1} > \theta_{j_2}$ if $j_1 < j_2$ and $j_1$ and $j_2$ lie in distinct intervals $[p_{k_1} + 1, p_{k_1} + q_{k_1}]$ and $[p_{k_2} + 1, p_{k_2} + q_{k_2}]$. (4.1)

Given a function $f$ of two variables, and functions $g_1$ and $g_2$ of one variable, write $\iint f g_1 g_2$ for $\iint f(u,v) g_1(u) g_2(v) \, du \, dv$. In this notation, and provided that

$$\sup_{u,v \in \mathcal{I}} |K(u,v)| < \infty \quad \text{and} \quad \int E(X^4) < \infty, \tag{4.2}$$

the random variables

$$M_{j_1 j_2} = n^{1/2} \iint (\widehat{K} - K) \psi_{j_1} \psi_{j_2}, \qquad j_1, j_2 \geq 1, \quad \text{and}$$
$$M_0 = n^{1/2} \sum_{j=1}^{\infty} (\hat{\theta}_j - \theta_j) = n^{1/2} \int (\widehat{K} - K)(u,u) \, du \tag{4.3}$$

are jointly, asymptotically distributed as multivariate normal variables $N_{j_1 j_2}$ and $\sum_j N_{jj}$, say, with zero means. In particular, if the principal components $\xi_j$ are as defined at (2.5), and if $Z_0 = \sum_j N_{jj}$, then

$$\text{cov}(N_{j_1 j_2}, N_{j_3 j_4}) = \text{cov}(\xi_{j_1} \xi_{j_2}, \xi_{j_3} \xi_{j_4}),$$

$$\operatorname{cov}(N_{j_1 j_2}, Z_0) = \sum_{j=1}^{\infty} \operatorname{cov}(\xi_{j_1} \xi_{j_2}, \xi_j^2). \quad (4.4)$$

Given $1 \leq k < \nu$, let $\mathcal{J}_k$ denote the set of indices $p_k + 1, \ldots, p_k + q_k$. Then $\theta_j = \theta^{(k)}$, say, not depending on $j$, for all $j \in \mathcal{J}_k$. Theorem 1 below asserts that the limiting joint distribution of the normalized eigenvalue differences $n^{1/2}(\hat{\theta}_j - \theta_j)$, and of $M_0$, is that of random variables $Z_j$, defined as follows: If $\mathcal{J}_k$ is a singleton, or equivalently, if there are no ties for the eigenvalue $\theta_{p_k+1}$, then we take $Z_{p_k+1} = N_{p_k+1,p_k+1}$. If $\mathcal{J}_k$ contains $q_k \geq 2$ elements, let $N^{(k)}$ denote the $q_k \times q_k$ matrix with $(j_1, j_2)$th component $N_{j_1 j_2}$, for $p_k + 1 \leq j \leq p_k + q_k$. Write $\Lambda_{p_k+1} \leq \cdots \leq \Lambda_{p_k+q_k}$ for the ordered eigenvalues of $N^{(k)}$. Since $k < \nu$, then, with probability 1, these eigenvalues are also distinct. Define $Z_{p_k+j} = \Lambda_{p_k+j}$ for $p_k + 1 \leq j \leq p_k + q_k$.

**Theorem 1.** *Assume* (4.1) *and* (4.2), *and that* $1 \leq k < \nu$. *Then the random variables* $n^{1/2}(\hat{\theta}_j - \theta_j)$, *for* $1 \leq j \leq p_k + q_k$, *and* $M_0$, *are jointly asymptotically distributed as* $Z_j$, *for* $1 \leq j \leq p_k + q_k$, *and* $Z_0$.

### 4.2. Properties of tie diagnostics

In Section 2.3 we suggested a particular member of a large class of diagnostics, the general class characterized by the property:

$$\text{decide that } \theta_{j_1} = \theta_{j_2} \quad \text{if and only if} \quad |\hat{\theta}_{j_1} - \hat{\theta}_{j_2}| \leq \hat{z}. \quad (4.5)$$

The value of the "critical point" $\hat{z}$ would generally be computed from data. If $\hat{z}$ were to satisfy the relations

$$0 \leq \hat{z} = o_p(1) \quad \text{and} \quad n^{-1/2} = o_p(\hat{z}), \quad (4.6)$$

as $n \to \infty$ and if the tie diagnostic were given by (4.5), then from Theorem 1 it would follow that the tie diagnostic asymptotically correctly identified a finite number of eigenvalue clusters, that is,

$$P(\hat{p}_\ell = p_\ell \text{ and } \hat{q}_\ell = q_\ell \text{ for } 1 \leq \ell \leq k) \to 1 \quad (4.7)$$

as $n \to \infty$ for $1 \leq k < \nu$. Any sequence $\hat{z}$ that decreased to zero more slowly than $n^{-1/2}$ would satisfy (4.6). Theorem 2, below, implies that the sequence $\hat{z} = 2\hat{z}_\beta$, introduced in Section 2.3 using a bootstrap argument, satisfies (4.6) if $\beta$ is permitted to decrease slowly to 0 with increasing sample size. In Section 3 we show that under (4.1) and (4.2) property (4.7) is sufficient for the consistency of TRB.

To state Theorem 2, let $W_{j_1 j_2}$ denote jointly normally distributed random variables with zero means and the same covariance structure as $\xi_{j_1} \xi_{j_2} - E(\xi_{j_1} \xi_{j_2})$. Put $\Theta_k^2 =$ $\sum_{j \in \mathcal{J}_k} (\hat{\theta}_j - \theta^{(k)})^2$, where $\theta^{(k)}$ denotes the common value of $\theta_j$ for $j \in \mathcal{J}_k$, and define

$$W_k^2 = \sum_{j_1 \in \mathcal{J}_k} \sum_{j_2 \in \mathcal{J}_k} W_{j_1 j_2}^2, \qquad W^2 = \sum_{j_1=1}^{\infty} \sum_{j_1=1}^{\infty} W_{j_1 j_2}^2.$$

Recall the definition of $\widehat{K}^\dagger$ from Section 2.3.

**Theorem 2.** *Assume (4.1) and (4.2), and that $1 \leq k < \nu$. Then, (a) $E(W) < \infty$, (b)*

$$n(\Theta_1^2, \ldots, \Theta_k^2, \|\widehat{K} - K\|^2) \to (W_1^2, \ldots, W_k^2, W^2), \tag{4.8}$$

*where the convergence is in joint distribution, and (c) conditional on the data $X_1, \ldots, X_n$, we have $n\|\widehat{K}^\dagger - \widehat{K}\|^2 \to W^2$ in distribution.*

Parts (b) and (c) of Theorem 2 imply that, as $n \to \infty$, the scaled bootstrap critical point $n^{1/2} \hat{z}_\beta$ converges in probability to the $(1-\beta)$-level critical point of $W$. Also, part (b) makes clear in an asymptotic sense the extent of conservatism of the bound

$$\sup_{1 \leq j \leq p_k + q_k} |\hat{\theta}_j - \theta_j| \leq \|\widehat{K} - K\|, \tag{4.9}$$

which forms the basis for the tie diagnostic suggested in Section 2.3. To appreciate this point, note that the square of the left-hand side of (4.9) is bounded above by $\max(\Theta_1^2, \ldots, \Theta_k^2)$, which in turn is bounded above by $\sum_{\ell \leq k} \Theta_\ell^2$, which, when multiplied by $n$, converges in distribution to $\sum_{\ell \leq k} \sum_{j_1 \in \mathcal{J}_\ell} \sum_{j_2 \in \mathcal{J}_\ell} W_{j_1 j_2}^2$. This triple series is, with probability 1, strictly less than $W^2$, which equals the limit, as $n \to \infty$, of $n\|\widehat{K} - K\|^2$.

### 4.3. Consistency of tie-respecting bootstrap

The following theorem implies consistency of TRB estimators of the joint distribution of any finite number of the $\hat{\theta}_j$'s, and of $\sum_j \theta_j$.

**Theorem 3.** *Assume (4.1) and (4.2) and that $1 \leq k < \nu$. If the tie diagnostic algorithm used in Section 2.6 to construct the bootstrap algorithm asymptotically correctly identifies the eigenvalue clusters up to the kth, that is, if it satisfies (4.7), then the joint distribution of $n^{1/2}(\hat{\theta}_j^* - \hat{\theta}_j)$ for $1 \leq j \leq p_k + q_k$ and $n^{1/2} \sum_j (\hat{\theta}_j^* - \hat{\theta}_j)$, conditional on the data $X_1, \ldots, X_n$, also converges to the joint distribution of $Z_1, \ldots, Z_{p_k+q_k}$ and $Z_0$.*

An immediate corollary is that simultaneous confidence regions for any finite number of $\theta_j$'s, and for $\sum_k \theta_k$ and the ratio $\sum_{k \leq j} \theta_k / \sum_k \theta_k$, have asymptotically correct coverage. Examples of such confidence regions include those at (2.14) and (2.15).

We have seen that confidence regions based on our TRB algorithm and the tie diagnostic in Section 2.3 have asymptotically correct coverage if $\beta$ decreases slowly to 0 with increasing sample size. However, the case of fixed $\beta$ is arguably both simpler and more easily interpretable. Its treatment makes explicit the cost of uncertainty in the decision made at the tie diagnostic stage. To appreciate why, note that Theorem 2 implies that the bootstrap critical point $\hat{z}_\beta$ is asymptotically conservative, in the sense that, for fixed $\beta$,

$$\lim_{n\to\infty} P\left(\sup_{j\geq 1} |\hat{\theta}_j - \theta_j| \leq \hat{z}_\beta\right) > 1 - \beta. \tag{4.10}$$

The strictness of the above inequality follows from the strictness of the bound

$$\sum_{\ell \leq k} \sum_{j_1 \in \mathcal{J}_\ell} \sum_{j_2 \in \mathcal{J}_\ell} W_{j_1 j_2}^2 < W^2,$$

discussed immediately after Theorem 2. From (4.10) it follows that, instead of (4.7),

$$\liminf_{n\to\infty} P(\hat{p}_\ell = p_\ell \text{ and } \hat{q}_\ell = q_\ell \text{ for } 1 \leq \ell \leq k) > 1 - \beta. \tag{4.11}$$

Bearing in mind that the limit property in Theorem 3 is conditional on the data, of which the quantities $\hat{p}_k$ and $\hat{q}_k$ are functions, it can be seen from (4.11) that a nominal $(1-\alpha)$-level confidence region based on our bootstrap algorithm, and using the tie diagnostic in Section 2.3 for fixed $\beta$, has asymptotic coverage strictly greater than $1-\alpha-\beta$. The strictness of this bound provides a degree of assurance of the robustness of our tie diagnostic, as we observed in Section 3.

### 4.4. Uniformity and perturbations of the model

In general the bootstrap is not good at producing distribution estimators and confidence regions that perform well in a uniform sense. See, for example, Hall and Jing (1995) and Romano (2004). Results of Eaton and Tyler (1991) and Dümben (1993) indicate that this is also a challenge for $m$-out-of-$n$ bootstrap methods. However, we are not aware of other methods that outperform the bootstrap in this respect, and the problem is probably inherent, rather than a shortcoming of the bootstrap *per se*.

Issues of this type arise if we consider perturbations of the covariance model as $n$ increases. To explore this problem, consider the case where we have a triangular array of functions, $X_{n,1}, X_{n,2}, \ldots, X_{n,n}$ for $n \geq 1$, with covariance kernel $K_n$ and eigenvalues $\theta_{n,1} \geq \theta_{n,2} \geq \cdots \geq 0$. Suppose that $K_n - K$ and $\theta_{nj} - \theta_j$ converge to zero at rate $\delta_1 = \delta_1(n)$, denoting a positive sequence decreasing to zero, and that $\hat{z}$, in (4.6), converges to zero at rate $\delta_2$. It will be assumed that $n^{-1/2}$ is of smaller order than both $\delta_1$ and $\delta_2$. If $\delta_1$ is of strictly larger order than $\delta_2$, that is, if $\delta_2 = o(\delta_1)$, then the probability that $|\hat{\theta}_{n,j_1} - \hat{\theta}_{n,j_2}| > \hat{z}$ converges to 1 as $n \to \infty$. Equivalently, with probability converging to 1 the tie diagnostic declares $\theta_{j_1}$ and $\theta_{j_2}$ to be unequal. This is the correct decision in the present case, since the rate $n^{-1/2}$ at which the eigenvalues $\theta_{n,j_1}$ and $\theta_{n,j_2}$ are estimated by $\hat{\theta}_{n,j_1}$ and $\hat{\theta}_{n,j_2}$, respectively, is of strictly smaller order than $\delta_1$, and so the limiting distributions of $\hat{\theta}_{n,j_1}$ and $\hat{\theta}_{n,j_2}$ are those that arise for unequal eigenvalues. However, if $\delta_1$ is of strictly smaller order than $\delta_2$, that is, if $\delta_1 = o(\delta_2)$, then the probability that $|\hat\theta_{n,j_1} - \hat\theta_{n,j_2}| \leq \hat z$ converges to 1 as $n \to \infty$, and so with probability converging to 1 the tie diagnostic declares $\theta_{j_1}$ and $\theta_{j_2}$ to be equal. This is the incorrect decision, and in this instance the TRB generally gives inconsistent estimators of the distributions of $\hat\theta_{n,j_1}$ and $\hat\theta_{n,j_2}$.

## 5. Technical arguments

### 5.1. Preliminary lemma

Assume that the kernels, $K_1$ and $K_2$, of two positive semi-definite Hilbert–Schmidt operators, also denoted by $K_1$ and $K_2$, are bounded and admit the spectral decompositions $K_k(u,v) = \sum_j \theta_{kj} \psi_{kj}(u) \psi_{kj}(v)$, where the $\theta_{kj}$'s are eigenvalues and the $\psi_{kj}$'s are the respective orthogonal eigenvectors. Suppose, too, that the eigenvalues are arranged in order of decreasing size, and that the eigenvalues of $K_1$ can be grouped into blocks of length $q_k$, where:

$q_k \geq 1$ and $p_k + q_k = p_{k+1}$ for each $k$, and $1 = p_1 + 1 \leq p_1 + q_1 < p_2 + 1 \leq p_2 + q_2 < p_3 + 1 \leq p_3 + q_3 < \cdots$, with the sequence of inequalities either continuing ad infinitum or ending with a value of $\nu$ for which $q_\nu = \infty$; $\theta_{1j_1} = \theta_{1j_2}$ (5.1) for $p_k + 1 \leq j_1, j_2 \leq p_k + q_k$ and each $k$; and $\theta_{1j_1} > \theta_{1j_2}$ if $j_1 < j_2$ and $j_1$ and $j_2$ lie in distinct intervals $[p_{k_1}+1, p_{k_1}+q_{k_1}]$ and $[p_{k_2}+1, p_{k_2}+q_{k_2}]$.

Put $\|\psi\|^2 = \int \psi^2$ and $\eta^2 = \|K_1 - K_2\|^2 = \iint (K_1 - K_2)^2$; for $r = 1$ or 2, define $\bar\psi_{2j} = \sum_{p+1 \leq \ell \leq p+q} \psi_{1\ell} \int \psi_{1\ell} \psi_{2j}$ and let $\bar\psi_{2j}$ denote the projection of $\psi_{2j}$ onto the space of functions spanned by $\psi_{1,p_k+1}, \ldots, \psi_{1,p_k+q_k}$.

Given $k \geq 1$, let $\mathcal{J}_k$ be the set of indices $j$ such that $p_k + 1 \leq j \leq p_k + q_k$. Let $B$ denote the $q_k \times q_k$ matrix with $(j_1, j_2)$th component $b_{j_1 j_2} = \eta^{-1} \iint (K_2 - K_1) \psi_{1j_1} \psi_{1j_2}$. Then $B$ is symmetric, and its eigenvalues, $\gamma_{p_k+1}, \ldots, \gamma_{p_k+q_k}$, say, are all real, with respective orthonormal eigenvectors $c_{p_k+1}, \ldots, c_{p_k+q_k}$ satisfying $Bc_j = \gamma_j c_j$ for $p_k + 1 \leq j \leq p_k + q_k$. The sum of the squares of the components of $B$ is bounded by 1, and so $\|c_j^\mathrm{T} B\| \leq 1$, implying that $|\gamma_j| \leq 1$ for each $j$.

Write $c_j = (c_{j,p_k+1}, \ldots, c_{j,p_k+q_k})^\mathrm{T}$. Let $e_k > 0$ denote a lower bound to the minimum spacing between adjacent $\gamma_j$'s; note that $e_k$ is assumed strictly positive. Observe that, by definition of $\mathcal{J}_k$, $\theta_{1j} = \theta^{(k)}$, say, not depending on $j$, for all $j \in \mathcal{J}_k$. Order the values of $\gamma_{p_k+1}, \ldots, \gamma_{p_k+q_k}$ as $\gamma_{(p_k+1)} \leq \cdots \leq \gamma_{(p_k+q_k)}$. Let $\pi$ denote the permutation of $p_k + 1, \ldots, p_k + q_k$ that takes $p_k + j$ to $\pi(p_k + j) = p_k + \ell_j$, defined such that $\gamma_{p_k+\ell_j} = \gamma_{(p_k+j)}$. Ties can be broken arbitrarily; it follows from (5.2), below, that if $e_k > 0$ is kept fixed, then ties do not occur if $\eta$ is sufficiently small. A proof of the lemma below is given in a longer version of this paper (Hall *et al.* (2007)).

**Lemma.** *If* (*5.1*) *holds, then for each* $k \geq 1$ *for which* $q_k$ *is finite there exist constants* $C_{k8}, C_{k9} > 0$, *depending only on the eigenvalue sequence* $\theta_{11}, \ldots, \theta_{1,p_k+q_k+1}$, *on* $e_k > 0$ and on $|\mathcal{I}|^2 \sup |K_1|$, *such that, for* $\eta \leq C_{k8}$ *and* $p_k + 1 \leq j \leq p_k + q_k$,

$$|\theta_{2j} - \theta^{(k)} - \eta \gamma_{\pi(j)}| \leq C_{k9} \eta^2. \tag{5.2}$$

## 5.2. Proofs of Theorems 1 and 3

Theorem 1 follows from the lemma on taking $K_1 = K$ and $K_2 = \widehat{K}$, whereas Theorem 3 follows on letting $K_1 = \widehat{K}$ and $K_2 = \widehat{K}^*$. In the context of Theorem 1, to treat the joint distribution of $n^{1/2}(\hat{\theta}_j - \theta_j)$ (for any finite number of values of $j$) and $M_0 = n^{1/2} \sum_j (\hat{\theta}_j - \theta_j)$, we note that the latter quantity has the simple and explicit representation in (4.3), which simplifies asymptotically to:

$$M_0 = \int \left\{ n^{-1/2} \sum_{i=1}^n X_i(u)^2 - K(u,u) \right\} du + o_p(1)$$

$$= \sum_{j=1}^\infty \left\{ n^{-1/2} \sum_{i=1}^n (\xi_{ij}^2 - \theta_j) \right\} + o_p(1).$$

Therefore, $M_{j_1 j_2}$ (for any finite collection of values of $(j_1, j_2)$) and $M_0$, both defined at (4.3), are jointly asymptotically distributed as $N_{j_1 j_2}$ and $Z_0$; see (4.4). Theorem 3 may be treated similarly.

## 5.3. Proof of Theorem 2

Assume, without loss of generality, that $\mathrm{E}(X) = 0$, and define $\xi_{ij} = \int X_i \psi_j$, $V_{ij_1 j_2} = \xi_{ij_1} \xi_{ij_2} - \delta_{j_1 j_2} \theta_{j_1}$, $U_{j_1 j_2} = n^{-1} \sum_i V_{ij_1 j_2}$ and

$$\widehat{L}(u,v) = \frac{1}{n} \sum_{i=1}^n \{X_i(u) X_i(v) - K(u,v)\}.$$

Then,

$$\|\widehat{L}\|^2 = \sum_{j_1=1}^\infty \sum_{j_1=1}^\infty U_{j_1 j_2}^2, \qquad |\|\widehat{K} - K\| - \|\widehat{L}\|| \leq \int \bar{X}^2 = \mathrm{O}_p(n^{-1}). \tag{5.3}$$

Let $\xi_j$ be as in (2.5). The second part of (4.2) implies that

$$\iint \mathrm{E}(XX - K)^2 = \sum_{j_1=1}^\infty \sum_{j_1=1}^\infty \mathrm{E}(\xi_{j_1} \xi_{j_2} - \delta_{j_1 j_2} \theta_{j_1})^2 < \infty. \tag{5.4}$$

Since $\mathrm{E}(W)$ equals the double series in (5.4), then (5.4) implies part (a) of Theorem 2.

Note that

$$\sum_{j_1=r_1+1}^{\infty}\sum_{j_2=r_2+1}^{\infty} n\mathrm{E}(U_{j_1j_2}^2) = \sum_{j_1=r_1+1}^{\infty}\sum_{j_2=r_2+1}^{\infty} \mathrm{E}(\xi_{j_1}\xi_{j_2} - \delta_{j_1j_2}\theta_{j_1})^2. \qquad (5.5)$$

Together, (5.4), (5.5) and the first part of (5.3) imply that for each $\varepsilon > 0$

$$\lim_{r_1,r_2\to\infty}\limsup_{n\to\infty} P\left(n\left|\|\widehat{L}\|^2 - \sum_{j_1=1}^{r_1}\sum_{j_2=1}^{r_2} U_{j_1j_2}^2\right| > \varepsilon\right) = 0. \qquad (5.6)$$

Define $U = (u_{j_1j_2})$ to be a $q_k \times q_k$ matrix, with $j_1, j_2 \in \mathcal{J}$ and $u_{j_1j_2} = n^{1/2}\iint \widehat{L}\psi_{j_1}\psi_{j_2} = n^{1/2}U_{j_1j_2}$. It can be shown from (5.2), on taking $K_1 = K$ and $K_2 = \widehat{K}$, that if $k < \nu$, then the sum of $n(\hat{\theta}_j - \theta^{(k)})^2$ over $j \in \mathcal{J}_k$ equals the sum of the eigenvalues of $U^2$, identical to the trace of $U^2$, plus a remainder $\mathrm{o}_p(1)$. Of course, the trace of $U^2$ equals $\sum_{j_1\in\mathcal{J}_k}\sum_{j_2\in\mathcal{J}_k} u_{j_1j_2}^2$. Therefore,

$$\sum_{j\in\mathcal{J}_k}(\hat{\theta}_j - \theta^{(k)})^2 = \sum_{j_1\in\mathcal{J}_\ell}\sum_{j_2\in\mathcal{J}_\ell} U_{j_1j_2}^2 + \mathrm{o}_p(n^{-1}). \qquad (5.7)$$

It can be proved using a conventional central limit theorem that, for each fixed $s_1, s_2 \geq 1$, the joint distribution of $n^{1/2}U_{j_1j_2}$, for $1 \leq j_\ell \leq s_\ell$ and $\ell = 1, 2$, converges weakly to the joint distribution of $W_{j_1j_2}$, for $j_1, j_2$ in the same range. This property, (5.6), (5.7) and both parts of (5.3) imply (4.8), and hence also part (b) of Theorem 2. Part (c) follows from a bootstrap version of these arguments, which is similar.

## Acknowledgements

Research of Byeong U. Park was supported by the Korea Science and Engineering Foundation (KOSEF) grant funded by the Korea government (MOST) (No. R01-2007-000-10143-0).